\documentclass[review]{elsarticle}
\usepackage{amsthm,amssymb,amsmath,amsfonts}
\usepackage{cases}
\usepackage{graphicx}
\usepackage{graphicx}
\usepackage{tikz-cd} 
\usepackage{lineno,hyperref}
\modulolinenumbers[5]
\theoremstyle{definition}
\newtheorem{definition}{Definition}[section]

\newtheorem{remark}{remark}[section]
\journal{ArXiv}

\begin{document}

\begin{frontmatter}

\title{A new method for compactification with the help of order topology and limit point}

\author{Kaveh Mohammadi\fnref{myfootnote}}
\address{Amirkabir University of Technology}
\fntext[myfootnote]{kmsmath@aut.ac.ir}

\author{Assad Rashidi\fnref{myfootnote}}
\address{University of Kurdistan}
\fntext[myfootnote]{a.rashidi@sci.uok.ac.ir}

\begin{abstract}
In this paper, we introduce a new method for compactification of a topological space by order topology and through ordinal numbers.  The idea behind our approach originates from the definition of a limit point, and then we try to find an intuition for this concept. Finally, we utilise the Homotopy concept for separation Axiom.
\end{abstract}


\end{frontmatter}

\section{Introduction}
In this paper, we introduce a new method for compactification of a topological space by order topology and through ordinal numbers.  The idea behind our way originates from the definition of a limit point, and then we try to find an intuition for this concept. Finally, we utilise the Homotopy concept for seperation Axiom. This paper is organised into four chapters. In the first chapter, we explain that the behaviour of limit points resembles a black hole in the branch of physics. In the second chapter, we introduce our method of compactification in four steps and with five standard shapes. In the third chapetr, with the help of homotopy theory, we explain the separation axiom for our approach. In the final chapter, we compare our method with other methods of compactification, such as one-point compactification and stone-Cech compactification.

\section{Chapter one}
\begin{definition}[limit point]

 Let A be a subset of a topological space W. A point w in W is a limit point of A if every neighbourhood of w contains at least one point  of A different than w.
\end{definition}
\begin{remark}
	A very prominent point in here is that when w doesn't have to be an element of A, in order to be a limit point of A and when a limit point exists of the set of A we have no problem.
	
	\begin{figure}[h!]
		\includegraphics[width=3cm,height=3cm,keepaspectratio]{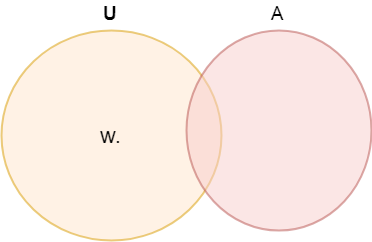}
		
	\end{figure}
When w belongs to A a serious question will be raised here " How is it possible to imagine every neighbourhood of w in the A doesn't have to contain w while w is a member of A.  Our idea for answering this question is that if we think of this point as a hole in the space then the problem of imagining a limit inside the set Of A will be resolved

\begin{figure}[h!]
	\includegraphics[width=3cm,height=3cm,keepaspectratio]{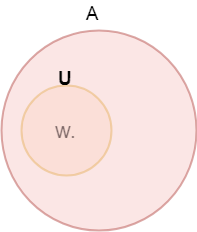}

\end{figure}
\end{remark}
\begin{definition}
	Let A be a subset of W a point $w\in \bar{W}$ is a limit (accumulation) point of A if every neighbourhood of $w$ contains at least one point of A different from $w$ itself.
\end{definition}

\begin{remark}
	Another concerning question is raised here is that in the second definition of limit point it is mentioned that if accumulation  point $w$  exists outside the set  for every neighbourhood of $w$ contains at least one point of A different from w itself but the issue in here is that we cannot imagine this point because for every neighbourhood of limit point w like V1 that
	$V_1\cap A\neq \emptyset$
	 we can find a neighbourhood like $V_{-\infty}$ such that $V_{-\infty} A=\emptyset$ because the existence of  distance between accumulation point w and subset A is continual. To address this issue, we can think of limit point as a black hole in the space, and this point attracts subset A toward itself hence we can imagine that for every neighbourhood of limit point $w$ like  $V_i,( i=1 \quad \text{to}  \quad i=-\infty)$we have  $V_i\cap A \neq \emptyset$
	 \begin{figure}[h!]
	 	\includegraphics[width=3cm,height=3cm,keepaspectratio]{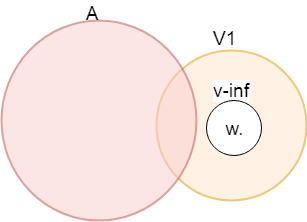}
	 	
	 \end{figure}
\end{remark}

\begin{remark}
	Another reason for this interpretation  of accumulation point is one-point compactification of the real line into a circle with one missing point infinity

	\begin{figure}[h!]
		\includegraphics[width=10cm,height=8cm,keepaspectratio]{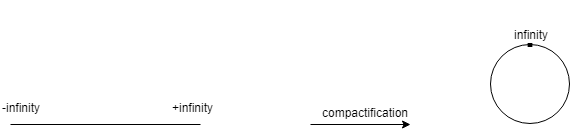}
		
	\end{figure}

	This illustration is a clear example for illumination our idea for thinking of accumulation point as a black hole in which the point at infinity attracts both ends of the real line toward itself, and then a circle is formed.
\end{remark}
\begin{remark}
	Based on the  connection  we found between black hole and limit-points we  have three following properties
	
	\begin{itemize}
	    \item [(a)] Limit-point is a hole in space.
	\item[(b)] Limit-points attracts sets in the space toward itself.
	\item[(c)] Limit points first draws each other, then bend their surrounding space
	\end{itemize}
\end{remark}

\section{Chapter 2}

I this chapter we are going to explain our compactification method in four steps.

\subsection{Step 1}

\begin{figure}[h!]
	\includegraphics[width=10cm,height=8cm,keepaspectratio]{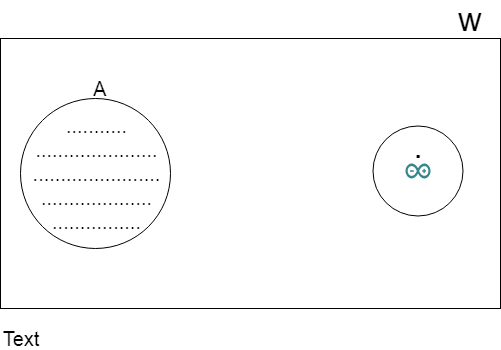}
	
\end{figure}
Let suppose that W is a compact  space and there is a subset of it which we called it A, and We want to embed this subset as dense subset of W in which $\infty$ is an accumulation point of A. We will show that $\bar{A}=W$. In this situation, we know that every point like w in W either belongs to A or is a limit point of A. The first condition is satisfied(If Space A doesn't contain all point of W other thank limit points Then we cannot compactify A) and for the second condition we must do the following steps to get close to the limit point infinity of the compact space W.

\subsection{Step 2}

\begin{figure}[h!]
	\includegraphics[width=14cm,height=12cm,keepaspectratio]{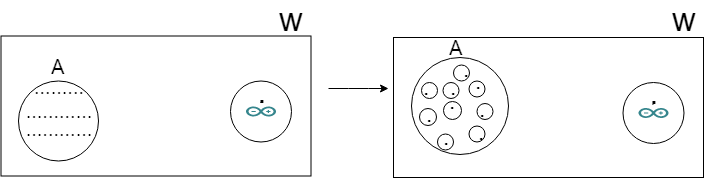}
	
\end{figure}
In this step we are going to crash the space A. In order to acheive this  purpose we know that we can put discrete or indiscrete topology on every set then we define an induction topology function:

\[C_i:A\rightarrow \{\tau_i\}\]
\[C_i(A)=\tau_i\]
such that $C_0(A)=\tau_0$(trivial topology) and $C_1(A)=\tau_1$(discrete topology), $i\in I=[0,1]$, ${\tau_i}$ (set of all topologies on the set A)
\begin{center}
	\begin{tikzcd}
	\tau_1 & \tau_0                                                                                                                                                & \tau_i \\
	\tau_i & A \arrow[r, "C_i"'] \arrow[u, "C_0"'] \arrow[ru, "C_i"'] \arrow[l, "C_i"'] \arrow[lu, "C_1"'] \arrow[ld, "C_i"'] \arrow[d, "c_i"'] \arrow[rd, "C_i"'] & \tau_i \\
	\tau_i & \tau_i                                                                                                                                                & \tau_i
	\end{tikzcd}
\end{center}

Remark: Function $C_i$ is continuous because at every moment in $I$ for A as an object we have one topology as another object then based of the definition of continuity in topological space, This function is continuous. Now, we can make a discrete space of A by continuous function $C_i$ in the process of converting topologies on the subset A from trivial to the discrete topology in which every point is isolated or as a clopen set.

\subsection{Step 3}
\begin{figure}[h!]
	\includegraphics[width=10cm,height=8cm,keepaspectratio]{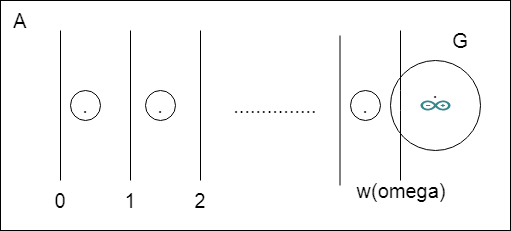}
	
\end{figure}
In this step want to show that $G\cap A\neq \emptyset$ to achieve this purpose we apply the ordinal topology with the help of ordinal numbers to space the  A so that it will become approximately  close to the accumulation point $\infty$.  By applying the map $\mathcal{O}$ on the point of discrete space A 
\[\mathcal{O}:[0,\omega)\rightarrow \tau^{A}_1\]
we obtain a ordinal-indexed sequence of point with discrete order topology over them in order to become close to the $\infty$.

\subsection{Step 4}
\begin{figure}[h!]
	\includegraphics[width=10cm,height=8cm,keepaspectratio]{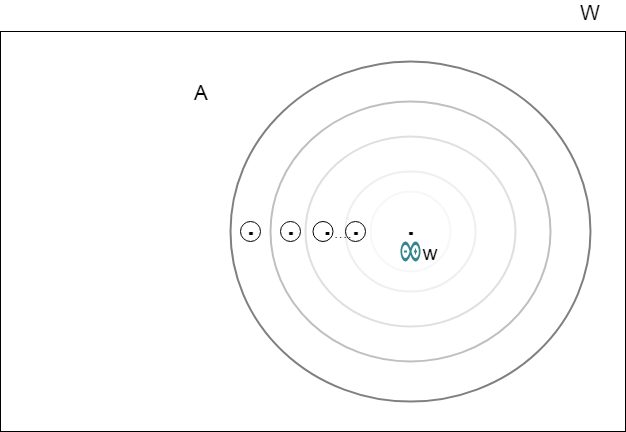}
	
\end{figure}

This is the final phase of the process of turning subset A into a dense subset of W. In regard with previous chapters limit pint $\infty$ acts like a black hole and it draws members of to itself and then ordinal space $[0,\omega)$ will turn into a compact and bounded space $[0,\infty_\omega]$ such that $\omega=\infty_\omega$.
\section{Chapter 3}

After introducing our method of compactification with the help of ordinal numbers and limit points acting as a black hole, we define this method of compactification for all separation axiom. So we set a homotopy between continuous  functions of $C_i$ such that $C_i:A\rightarrow {\tau_i}$ and then we have 
\[H:A\times I\rightarrow {\tau_i}\]
such that 
$H(A,0)=\tau_0$ (trivial topology)
$H(A,1)=\tau_1$ (indiscrete topology)

\begin{remark} Structure I:
	The set of all topologies on set A together with partially ordering relation $\subseteq$ forms a complete lattice and also bounded lattice that we can consider trivial topology as the least element and discrete topology as greatest element this algebraic structure of topology have the form $(I,\wedge,\vee,0,,1)$ such that 
	O is a trivial topology, and 1 is discrete topology.
	
	Hence every discrete topological space satisfies each of the separation axioms, and every indiscrete topological doesn't fill any of the separation axioms then topologies on every set have a direct relation with the separation axiom. We will make a connection between all of the topology on the subset A through homotopy, and then we explain all separation axion in regard with our method of compactification.
	
\end{remark}
We divide the process of separation into six tens interval in I and also this process can be divided into more interval for other separation axioms.

\begin{numcases}{H(A,i)=}
\tau_0\text{(trivial topology, have no separation axiom)} i=0 \\
\tau_i \text{(some topologies satiesfies $T_0$ axiom )} 0<i\leqslant \frac{1}{6} \\
\tau_i \text{(some topologies satiesfies $T_1$ axiom )} \frac{1}{6}<i\leqslant \frac{2}{6} \\
\tau_i\text{(some topologies satiesfies $T_2$ axiom )} \frac{2}{6}<i\leqslant \frac{3}{6} \\
\tau_i\text{(some topologies satiesfies $T_3$ axiom )} \frac{3}{6}<i\leqslant \frac{4}{6} \\
\tau_i\text{(some topologies satiesfies $T_4$ axiom )} \frac{4}{6}<i\leqslant \frac{5}{6} \\
\tau_i\text{(discrete topoloy satisfies $T_0,...,T_4$ )} \frac{5}{6}<i\leqslant 1 
\end{numcases} 
In consequence, based on the last step $i=1$ we have a discrete topology that satisfies all of the separation axioms so that we can use this method for all spaces with every separation axioms.

\section{Chapter 4}

In this chapter, we compare our method of compactification with other well-known methods of compactification.

\subsection{One-point compactification}
When $\omega =\infty_\omega$ the space A based on $[0,\omega)$ is just N with the usual order topology after using order topology compactification, we have $[0,\infty_\omega]$  which is one-point compactification of N.

\subsection{Stone-Cech compactification}

\begin{center}
\begin{tikzcd}
{[0,\omega)} \arrow[r, hook] \arrow[rd, "{c([0,\omega))}"'] & {[0,\infty_\omega]} \arrow[d, "{C([0,\infty_\omega])}"] \\
& R                                                    
\end{tikzcd}
\end{center}

$C([0,\omega))$: Real valued continuous functions from topological space $[0,\omega)$ to Real line.

$C([0,\infty_\omega])$: Real Valued continuous function from compact Hausdorf space $[0,\infty_\omega]$ to  the real line.

At $i=0$  (trivial topology) we have no separation axiom then any continuous function from $[0,\omega)$ to R is eventually constant,  the Stone-Cech compactification of $\omega$ is $[0,\infty_\omega]$.

\section{Refrence}


\begin{thebibliography}{}
\expandafter\ifx\csname url\endcsname\relax
  \def\url#1{\texttt{#1}}\fi
\expandafter\ifx\csname urlprefix\endcsname\relax\def\urlprefix{URL }\fi
\expandafter\ifx\csname href\endcsname\relax
  \def\href#1#2{#2} \def\path#1{#1}\fi

\end{thebibliography}


\medskip

\begin{thebibliography}{9}
	\bibitem{} 
 Allen Hatcher
	\textit{Algebraic Topology}. 
Cambridge Univ Pr; 1 edition (September 1, 2005)

\bibitem{l.gillmanm.jerison1960}
James Munkres. 
\textit{Topology; a first course}.
Pearson College Div (June 1, 1974).

\bibitem{sashokalajdzievski2015}
Sasho Kalajdzievski. 
\textit{An Illustrated Introduction to Topology and Homotopy}.
Chapman and Hall/CRC 2015.

\bibitem{da}
Assad Rashidi, Kaveh Mohammadi
\textit{Triple extension of Tietze theorem and Baer criterionn}.
arXiv preprint arXiv:1906.02784, 2019.

\bibitem{johndaun1994}
B. K. Lahiri 
\textit{A First Course in Algebraic Topology }.
Alpha Science International, Ltd; 1 edition (August 1, 2000).




\bibitem{johndaun1994}
B. K. Lahiri 
\textit{A First Course in Algebraic Topology }.
Alpha Science International, Ltd; 1 edition (August 1, 2000).

\bibitem{willardstephen2004}
Willard, Stephen. 
\textit{General Topology}.
Dover Publications 2004.







\end{thebibliography}
\end{document}